\documentclass[12pt]{amsart}

\newtheorem{theorem}{Theorem}
\newtheorem{lemma}{Lemma}
\newtheorem{proposition}{Proposition}
\newtheorem{corollary}{Corollary}
\newtheorem*{thm}{Theorem}
\newtheorem*{prop}{Proposition}
\theoremstyle{definition}
\newtheorem{definition}{Definition}
\theoremstyle{remark}
\newtheorem{remark}{Remark}

\newtheorem{example}{Example}

\newcommand{\ltwo}{L^2({\mathbb R})}
\newcommand{\ltwod}{L^2({\mathbb R}^d)}

\newcommand{\Cal}{\mathcal}
\newcommand{\Bb}{\mathbb}
\newcommand{\THG}{\Theta_{\Bb{H}, \Bb{G}}}

\begin{document}


\title{Orthogonal Frames of Translates}
\author{Eric Weber}
\address{Department of Mathematics, Iowa State University, 400 Carver Hall, Ames, IA 50011} 
\email{esweber@iastate.edu}
\thanks{This research was supported in part by NSF grant DMS-0308634 \\ Submitted to Appl. Comp. Harmonic Anal.}
\date{\today}
\begin{abstract}
Two Bessel sequences are orthogonal if the composition of the synthesis operator of one sequence with the analysis operator of the other sequence is the 0 operator.  We characterize when two Bessel sequences are orthogonal when the Bessel sequences have the form of translates of a finite number of functions in $\ltwod$.  The characterizations are applied to Bessel sequences which have an affine structure, and a quasi-affine structure.  These also lead to characterizations of superframes.  Moreover, we characterize perfect reconstruction, i.e. duality, of subspace frames for translation invariant (bandlimited) subspaces of $\ltwod$.
\end{abstract}
\maketitle



\section*{Introduction}

Frames for (separable) Hilbert spaces were introduced by Duffin and Schaeffer \cite{DS52a} in their work on non-harmonic Fourier series.  Later, Daubechies, Grossmann, and Meyer revived the study of frames in \cite{DGM86a}, and since then, frames have become the focus of active research, both in theory and in applications, such as signal processing.  Every frame (or Bessel sequence) determines an analysis operator, the range of which is important for a number of applications.  Information about this range is partially revealed by considering the composition of analysis and synthesis operators for different frames.  We view this composition as a sum of rank one tensors.  The present paper considers frames and Bessel sequences in $\ltwod$ which arise from translations of generating functions, such as in wavelet and Gabor frame theory.  The goal is to determine when the infinite sum of rank one tensors involving these translations is actually the $0$ operator.  See the subsection entitled "Motivation" below.

\subsection{Definitions}

Let $H$ be a separable Hilbert space and $\Bb{J}$ a countable index set.  A sequence $\Bb{X} := \{x_j\}_{j \in \Bb{J}}$ is a \emph{frame} if there exist positive real numbers $C_1$, $C_2$ such that for all $v \in H$,
\begin{equation} \label{E:frame}
C_1 \|v\|^2 \leq \sum_{j \in \Bb{J}} |\langle v, x_j \rangle|^2 \leq C_2 \|v\|^2.
\end{equation}
If $\Bb{X}$ satisfies only the second inequality, (i.e.~only $C_1 = 0$ satisfies the first inequality), then $\Bb{X}$ is called a \emph{Bessel sequence}.  Given $\Bb{X}$ which is Bessel, define the analysis operator
\[ \Theta_{\Bb{X}}: H \to l^2(\Bb{J}) : v \mapsto (\langle v, x_j \rangle)_{j};\]
and the synthesis operator
\[ \Theta_{\Bb{X}}^{*}: l^2(\Bb{J}) \to H : (c_j)_{j} \mapsto \sum_{j \in \Bb{J}} c_j x_j.\]
The analysis operator is well-defined and bounded by the frame inequality (\ref{E:frame}).  Additionally, the sum $\sum_{j} c_j x_j$ converges (see \cite{DS52a}), and so the synthesis operator is also well-defined and bounded, and a simple computation shows that it is in fact the adjoint operator of the analysis operator.

Given two Bessel sequences $\Bb{X}$ and $\Bb{Y} := \{ y_j\}_{j \in \Bb{J}}$, define the operator
\[ \Theta_{\Bb{Y}}^{*} \Theta_{\Bb{X}}: H \to H : v \mapsto \sum_{j \in \Bb{J}} \langle v, x_j \rangle y_j;\]
this operator is sometimes called a "Mixed Dual Grammian".  Note that it is a (convergent) sum of rank one tensors, as described above.  Typically in frame theory, one wants the above operator to be the identity; if this is the case, then the Bessel sequences $\Bb{X}$ and $\Bb{Y}$ are actually frames and are called \emph{dual frames}.  Our motivation here is for the operator to be the $0$ operator.

\begin{definition}
Suppose $\Bb{X}$ and $\Bb{Y}$ are Bessel sequences in $H$.  If 
\[ \Theta_{\Bb{Y}}^{*} \Theta_{\Bb{X}} := \sum_{j \in \Bb{J}} \langle \cdot , x_j \rangle y_j = 0,\]
the Bessel sequences are said to be \emph{orthogonal}.
\end{definition}
This idea has been studied by Han and Larson (\cite{HL00a}), where the Bessel sequences were assumed to be frames and were called strongly disjoint, and also by Balan, et. al. in \cite{B99a} and \cite{BL02a} for the Gabor (Weyl-Heisenberg) frame case.

Orthogonality also arises in the case of $M$-subspace frames (called outer frames in \cite{ACM03a}).  Here, we consider frames for subspaces in a Hilbert space where the elements of the frame are not necessarily elements of the subspace.
\begin{definition}
Suppose $M \subset H$ is a closed subspace, and let $\{x_j\} \subset H$.  If there exists constants $0 < C_1, C_2 < \infty$ such that for all $v \in M$,
\[ C_1 \| v \|^2 \leq \sum_{j \in \Bb{J}} | \langle v , x_j \rangle |^2 \leq \| v \|^2 ,\]
then $\{x_j\}$ is an $M$-subspace frame.  If $\{x_j\}$ and $\{y_j\}$ are Bessel sequences and for every $v \in M$,
\[ v = \sum_{j \in \Bb{J}} \langle v, x_j \rangle y_j,\]
then $\{y_j\}$ is an $M$-subspace dual to $\{x_j\}$.
\end{definition}

\begin{remark}
It is possible for $\{y_j\}$ to be $M$-subspace dual frame for $\{x_j\}$ and $\{x_j\}$ is NOT an $M$ subspace dual for $\{y_j\}$. See Example \ref{E:bidual} in section \ref{S:DFTIS} below.
\end{remark}

\begin{definition}
A Bessel sequence $\Bb{X} \subset H$ is a Plancherel frame for $M$ if for all $v \in M$, $v = \sum_{j \in \Bb{J}} \langle v, x_j \rangle x_j$.
\end{definition}

See also \cite{LO01a} for alternative duals.

\noindent
{\bf Notation.}
For the purposes of this paper, we will define the Fourier transform for $f \in L^1(\Bb{R}^d) \cap \ltwod$ to be
\[ \hat{f}(\xi) = \int f(x) e^{-2\pi i x \cdot \xi} dx. \]
Define the dense subspace $\Cal{D} \subset \ltwod$ to be
\[ \Cal{D} := \{ f \in \ltwod : \hat{f} \in L^{\infty}(\Bb{R}^d); \ supp(\hat{f}) \text{ is compact and bounded away from $0$} \}. \]
If $P \in B(H)$ is an orthogonal projection, let $P^{\perp}$ be the orthogonal projection such that $P + P^{\perp} = I$, the identity.  If $\Cal{A} \subset B(H)$, $\Cal{A}'$ denotes the commutant of $\Cal{A}$, that is
\[ \Cal{A}' = \{ B \in B(H) : AB = BA \ \forall A \in \Cal{A} \}.\]
Note that if $\Cal{A}$ is a self-adjoint collection of operators, then $\Cal{A}'$ is a von Neumann algebra.

If $C$ is an invertible real matrix, let $C' = C^{*-1}$, where $C^{*}$ is the transpose.

Finally, for $\alpha \in \Bb{R}^d$, let $T_{\alpha}$ denote the unitary translation operator
\[ T_{\alpha} : \ltwod \to \ltwod : f(\cdot ) \mapsto f(\cdot - \alpha).\]

\subsection{Motivation}
In both theory and applications it is desirable to know the range of the analysis operator for a given frame.  Consequently, it is desirable to know the orthogonal complement of the range.  This can be determined by considering which frames (and Bessel sequences) have orthogonal ranges.  We list here a few examples:
\begin{enumerate}
\item Duality:  In some applications, one wishes to know many duals to the fixed frame.  Let $\{x_j\}$ be a frame.  Suppose $\{y_j\}$ is a dual frame for $\{x_j\}$; hence $\Theta_{Y}^{*} \Theta_{X} = I$.  If $Z := \{z_j\}$ is Bessel and orthogonal to $\{x_j\}$, then $\{y_j + z_j\} =: Y + Z$ is also a dual to $\{x_j\}$:
\[ \Theta_{Y + Z}^{*} \Theta_{X} = \Theta_{Y}^{*} \Theta_{X} + \Theta_{Z}^{*} \Theta_{X} = I. \]
Conversely, if $\{w_j\}$ is dual to $\{x_j\}$, then $w_j = y_j + z_j$ for some orthogonal Bessel sequence $\{z_j\}$.  Hence, the orthogonal sequences parametrize all duals to a fixed frame.
\item Multiple Access Communications:  Suppose $\{x_j\} \subset H$ and $\{y_j\} \subset K$ are both Parseval frames and are orthogonal to each other.  Then for any $v \in H$ and $w \in K$, we have
\[ v = \sum (\langle v, x_j \rangle + \langle w, y_j \rangle) x_j \text{ and } w = \sum (\langle v, x_j \rangle + \langle w, y_j \rangle) y_j. \]
In other words, the frames can be used to encode two signals $v$ and $w$, which can then be sent over a single communications channel.  See \cite{Bal00a, BDV00a}.
\item Superframes: Superframes are frames of the form $\{ x_j \oplus y_j\} \subset H \oplus K$.  These are related to multiple access communications \cite{Bal00a}.
\item Perfect reconstruction in subspaces:  In some applications, notably sampling theory, frames for subspaces are used in which the frame elements are not actually in the subspace.  For example, when oversampling the bandlimited functions in the Paley-Wiener space, instead of reconstructing the function with the sinc function, which decays poorly, one can use a function $\phi$ such that $\hat{\phi}$ is smooth and is identically 1 on $[-1/2, 1/2]$ and decays sufficiently fast outside that band:  
\[ f(x) = \sum_{n} f(an) \phi(x - an). \]
This is only possible when the samples are faster than the Nyquist rate.  Moreover, the functions $\phi(x - an)$ are not in the Paley-Wiener space.  This perfect reconstruction is because of orthogonality of certain Bessel sequences (see section \ref{S:DFTIS}).  For similar results in sampling theory see \cite{A02a, Web02a}.
\end{enumerate}

\subsection{Main Results}
Here we will state a few representatives of the main results in the paper.  The main results center around the orthogonality of wavelet frames, the duality of wavelet frames, the characterization of Parseval superwavelets, and perfect reconstruction in subspaces.

\begin{thm}
Suppose $A$ is an expansive integral matrix and the affine systems generated by $\Psi = \{\psi_1, \dots, \psi_r\}$ and $\Phi = \{\phi_1, \dots, \phi_r\}$ with respect to the dilation matrix $A$ are both Bessel sequences.  Then they are orthogonal if and only if for all $q \in \Bb{Z}^d \setminus A^{*}\Bb{Z}^d$,
\[\sum_{i=1}^{r} \sum_{j \geq 0} \overline{\hat{\psi}_i(A^{*j} \xi)} \hat{\phi}_i(A^{*j}(\xi + q)) = 0 \ a.e. \ \xi,\]
and 
\[ \sum_{i=1}^{r} \sum_{j \in \Bb{Z}} \overline{\hat{\psi}_i(A^{*j} \xi)} \hat{\phi}_i(A^{*j}\xi) = 0 \ a.e. \ \xi. \]
Moreover, the corresponding quasi-affine sequences are orthogonal if and only if the same two equations hold.
\end{thm}

\begin{prop}
If $A$ and $B$ are any dilation matrices and the affine systems generated by $\Psi = \{\psi_1, \dots, \psi_r\}$ and $\Phi = \{\phi_1, \dots, \phi_r\}$ with respect to the dilation matrices $A$ and $B$, respectively are dual, then $A=B$.
\end{prop}

\begin{thm}
Suppose $A$ is an expansive integral matrix and the affine systems generated by $\psi_i$ with respect to the dilation matrix $A$ are Bessel sequences for $i = 1,\dots, r$.  The superwavelet generated by $\psi_1 \oplus \cdots \oplus \psi_r$ is a Parseval frame if and only if
\begin{enumerate}
\item $\sum_{n \in \Bb{Z}} \hat{\psi}_i(A^{*n} \xi) \overline{\hat{\psi}_j(A^{*n} \xi)} = \delta_{i,j} \ a.e \ \xi$ for $i,j=1,\dots,r$, and
\item $\sum_{n = 0}^{\infty} \hat{\psi}_i(A^{*n} \xi) \overline{\hat{\psi}_j(A^{*n} (\xi + k))} = 0 \ a.e \ \xi$ for $k \in \Bb{Z}^d \setminus A^{*}\Bb{Z}^d$ and $i,j=1,\dots,r$.
\end{enumerate}
\end{thm}

For perfect reconstruction in subspaces, see the following subsection and also Section \ref{S:DFTIS}.

\subsection{Preliminary Results}
For the purposes of this subsection, let $\Bb{X} = \{x_j\}_{j \in \Bb{J}}$ and
$\Bb{Y} = \{y_j\}_{j \in \Bb{J}}$ be sequences in $H$.

\begin{lemma} \label{L:adjoint}
Let $\Bb{X}$ and $\Bb{Y}$ be Bessel sequences, and let $\Theta = \sum_{j \in \Bb{J}} \langle \cdot , x_j \rangle y_j$.  Then $\Theta^{*} = \sum_{j \in \Bb{J}} \langle \cdot , y_j\rangle x_j$.
\end{lemma}
\begin{proof}
Let $v,w \in H$; since $\sum_{j \in \Bb{J}} \langle v ,x_j \rangle y_j$ converges in $H$, we have:
\[ \langle \Theta v, w \rangle = \langle \sum_{j \in \Bb{J}} \langle v, x_j \rangle y_j , w \rangle = \sum_{j \in \Bb{J}} \langle v, x_j \rangle \langle y_j, w \rangle = \sum_{j \in \Bb{J}} \overline{\langle w, y_j \rangle} \langle v, x_j \rangle = \langle v , \sum_{j \in \Bb{J}} \langle w, y_j \rangle x_j \rangle.\]
\end{proof}

\begin{lemma} \label{L:P_theta}
If $\Bb{X}$ and $\Bb{Y}$ are Bessel and $P$ is an orthogonal projection, then $ \Theta_{Y}^{*} \Theta_{X} \in \{P\}'$ if and only if $\sum_{j \in \Bb{J}} \langle \cdot , P x_j \rangle P^{\perp} y_j = 0$ and $\sum_{j \in \Bb{J}} \langle \cdot , P^{\perp} x_j \rangle P y_j = 0$.
\end{lemma}
\begin{proof}
Write
\begin{align*}
\sum_{j \in \Bb{J}} \langle \cdot , x_j \rangle y_j &= \sum_{j \in \Bb{J}} \langle \cdot , P x_j \rangle P y_j + \sum_{j \in \Bb{J}} \langle \cdot , P x_j \rangle P^{\perp} y_j + \sum_{j \in \Bb{J}} \langle \cdot , P^{\perp} x_j \rangle P y_j + \sum_{j \in \Bb{J}} \langle \cdot , P^{\perp} x_j \rangle P^{\perp} y_j \\
:&= A + B + C + D.
\end{align*}
Clearly, we have the following:
\[ PA = AP, \qquad PD = 0 = DP, \qquad BP = B, \qquad PB = 0, \qquad PC = C, \qquad CP = 0.\]
Therefore, since the range of $B$ is in $PH$ and the range of $C$ is in $P^{\perp}H$, 
\[ PA + PB + PC + PD = C \text{ and } AP + BP + CP + DP = B \]
are equal if and only if $B = C = 0$.
\end{proof}

\begin{lemma} \label{L:PXY}
Suppose $\Bb{X} := \{x_j:j \in \Bb{J}\}$ is a Bessel sequence in $H$ and let $P \in B(H)$ be the orthogonal projection onto the closed subspace $M \subset H$.  The collections $\{Px_j: j \in \Bb{J} \}$ and $\{P^{\perp} x_j: j \in \Bb{J} \}$ are orthogonal, i.e.~$\sum_{j \in \Bb{J}} \langle \cdot , Px_j \rangle P^{\perp}x_j = 0$ if and only if $\Theta_{\Bb{X}}^{*} \Theta_{\Bb{X}} \in \{P\}'$.
\end{lemma}
\begin{proof}
Apply Lemma \ref{L:P_theta} to $\{x_j\}$ and $\{y_j \} = \{x_j\}$.  Note that by Lemma \ref{L:adjoint} the adjoint operator of $\sum_{j \in \Bb{J}} \langle \cdot , P x_j \rangle P^{\perp} x_j$ is $\sum_{j \in \Bb{J}} \langle \cdot , P^{\perp}x_j \rangle P x_j$.
\end{proof}

\begin{lemma} \label{L:MSDF}
Suppose $\{x_j\}$ and $\{y_j\}$ are Bessel sequences; $\{y_j\}$ is an $M$-subspace dual frame for $\{x_j\}$ if and only if for every $v \in M$,
\begin{enumerate}
\item $v = \sum_{j \in \Bb{J}} \langle v , P x_j \rangle P y_j$, and
\item $0 = \sum_{j \in \Bb{J}} \langle v, P x_j \rangle P^{\perp}y_j$.
\end{enumerate}
\end{lemma}
\begin{proof}
Let $v \in M$ and consider
\[ \sum_{j \in \Bb{J}} \langle v , x_j \rangle y_j = \sum_{j \in \Bb{J}} \langle v , P x_j \rangle P y_j + \sum_{j \in \Bb{J}} \langle v , P x_j \rangle P^{\perp} y_j = v \]
if items 1.~and 2.~hold.

Conversely, suppose $\{y_j\}$ is an $M$-subspace dual frame for $\{x_j\}$.  Then for all $v \in M$,
\[ Pv = P \sum_{j \in \Bb{J}} \langle v, x_j \rangle y_j = \sum_{j \in \Bb{J}} \langle v, P x_j \rangle P y_j\]
and
\[ 0 = P^{\perp} v = P^{\perp} \sum_{j \in \Bb{J}} \langle v, x_j \rangle y_j = \sum_{j \in \Bb{J}} \langle v, P x_j \rangle P^{\perp}y_j. \]
\end{proof}

\begin{remark}
We remark again that it is possible for $\{y_j\}$ to be $M$-subspace dual frame for $\{x_j\}$ and $\{x_j\}$ is NOT an $M$ subspace dual for $\{y_j\}$. See Example \ref{E:bidual} in section \ref{S:DFTIS} below.  Note also that item 2.~above is equivalent to $0 = \sum_{j \in \Bb{J}} \langle v, P x_j \rangle P^{\perp}y_j$ for all $v \in H$.
\end{remark}

\begin{lemma}
Let $M \subset H$ be a closed subspace, let $P_M$ be the orthogonal projection onto $M$, and let $\Bb{X} \subset H$ be a Bessel sequence.  The following are equivalent:
\begin{enumerate}
\item $\Bb{X}$ is a Plancherel frame for $M$;
\item for all $v \in M$,
\begin{enumerate}
\item $\|v\|^2 = \sum_{j \in \Bb{J}} | \langle v, x_j \rangle |^2$;
\item $\sum_{j \in \Bb{J}} \langle v , x_j \rangle P_{M}^{\perp} x_j = 0$;
\end{enumerate}
\end{enumerate}
The following implies both 1 and 2:
\begin{enumerate}
\item[3.] for all $v \in M$,
\begin{enumerate}
\item $\|v\|^2 = \sum_{j \in \Bb{J}} | \langle v, x_j \rangle |^2$;
\item $\Theta_{\Bb{X}}^{*} \Theta_{\Bb{X}} \in \{P\}'$.
\end{enumerate}
\end{enumerate}
\end{lemma}
\begin{proof}
Suppose $\{x_j\}$ is a Plancherel frame for $M$.  Then clearly, the sequence $\{P_M x_j\}$ is a Parseval frame for $M$; whence it follows that
\[ \sum_{j \in \Bb{J}} | \langle v , x_j \rangle |^2 = \sum_{j \in \Bb{J}} | \langle v , P_M x_j \rangle |^2 = \| v \|^2. \]
Moreover, we have
\[ v = \sum_{j \in \Bb{J}} \langle v , x_j \rangle x_j = \sum_{j \in \Bb{J}} \langle v , x_j \rangle P_M x_j + \sum_{j \in \Bb{J}} \langle v , x_j \rangle P_M^{\perp}x_j.\]
Since $\sum_{j \in \Bb{J}} \langle v , x_j \rangle P_M^{\perp}x_j \in M^{\perp}$, it must be $0$. 

Conversely, if $\|v\|^2 = \sum_{j \in \Bb{J}} | \langle v, x_j \rangle |^2$, then for all $v \in M$, 
\begin{align*}
v &= \sum_{j \in \Bb{J}} \langle v , P_M x_j \rangle P_M x_j \\
&= \sum_{j \in \Bb{J}} \langle v , x_j \rangle P_M x_j + \sum_{j \in \Bb{J}} \langle v , x_j \rangle P_M^{\perp} x_j \\
&= \sum_{j \in \Bb{J}} \langle v , x_j \rangle x_j.
\end{align*} 

Finally, by Lemma \ref{L:P_theta}, condition 3(b) implies condition 2(b), whence condition 3 implies condition 2.
\end{proof}

Note that condition 2(b) is equivalent to $\sum_{j \in \Bb{J}} \langle \cdot , P_Mx_j \rangle P_M^{\perp} x_j = 0$.

\section{General Translation Systems}

As in \cite{HLW02a}, let $\Cal{P}$ be a countable index set, let $C_p$ be a $d \times d$ invertible matrix for each $p \in \Cal{P}$, and define the following:
\[ \Lambda = \cup_{p \in \Cal{P}} C_p^{\prime} \Bb{Z}^d \] 
and for $\alpha \in \Lambda$,
\[ \Cal{P}_{\alpha} = \{p \in \Cal{P}: C_p^{*} \alpha \in \Bb{Z}^d\}. \]
Note that if $\alpha = C_{p_0}^{\prime} z$ for some $z \in \Bb{Z}^d \setminus \{0\}$, then $p_0 \in \Cal{P}_{\alpha}$; if $\alpha = 0$, then $\Cal{P}_{\alpha} = \Cal{P}$.  Let $\{g_p:p \in \Cal{P}\} \subset \ltwod$.  The collection $\{T_{C_p k} g_p: p \in \Cal{P}, \ k \in \Bb{Z}^d\}$ satisfies the Bessel condition if there exists a constant $M < \infty$ such that for all $f \in \ltwod$,
\[ \sum_{p \in \Cal{P}} \sum_{k \in \Bb{Z}^d} | \langle f, T_{C_p k} g_p \rangle |^2 \leq M \|f \|^2. \]
The collection $\{T_{C_p k} g_p: p \in \Cal{P}, \ k \in \Bb{Z}^d\}$ satisfies the local integrability condition \cite{HLW02a} if for every $f \in \Cal{D}$,
\[ L(f) := \sum_{p \in \Cal{P}} \sum_{k \in \Bb{Z}^d} \int_{supp \hat{f}} |\hat{f} (\xi + C_{p}^{\prime} k) |^2 |\det C_p|^{-1} |\hat{g}_{p}(\xi)|^2 d \xi < \infty. \]

\begin{theorem} \label{T:g_h_alpha}
Suppose $\{T_{C_p k} g_p\}$ and $\{T_{C_p k} h_p\}$ satisfy the Bessel condition and the local integrability condition.  The operator
\[ \Theta := \theta_{g}^{*} \theta_{h} = \sum_{p \in \Cal{P}} \sum_{k \in \Bb{Z}^d} \langle \cdot , T_{C_p k} h_p \rangle T_{C_p k} g_p \]
is in the von Neumann algebra $\{T_{\beta}: \beta \in \Bb{R}^d\}'$ if and only if for all $\alpha \in \Lambda \setminus \{0\}$,
\begin{equation} \label{E:g_h_alpha}
 \sum_{p \in \Cal{P}_{\alpha}} |\det C_p|^{-1} \overline{\hat{h}_p(\xi)} \hat{g}_{p}(\xi + \alpha) = 0 \ a.e. \ \xi.\end{equation}
In this case, $\Theta$ is a Fourier multiplier whose symbol is
\[ s(\xi) = \sum_{p \in \Cal{P}} |\det C_p|^{-1} \overline{\hat{h}_{p}(\xi)} \hat{g}_{p}(\xi). \]
\end{theorem}
\begin{proof}
For $f \in \Cal{D}$, define the continuous function
\[ w_f(x) = \langle \Theta T_x f , T_x f \rangle. \]
If $\Theta$ commutes with all $T_{\beta}$ for $\beta \in \Bb{R}^d$, then clearly $w_f(x)$ is constant for all $f \in \Cal{D}$.  Conversely, if $w_f(x)$ is constant for all $f \in \Cal{D}$, then $\langle T_{-x} \Theta T_x f, f \rangle = \langle \Theta f, f \rangle$, whence by the polarization identity, $T_{-x} \Theta T_x = \Theta$, and thus $\Theta T_x = T_x \Theta$.

By \cite[Proposition 2.4]{HLW02a}, $w_f(x)$ coincides pointwise with the almost periodic function 
\[\sum_{\alpha \in \Lambda} \hat{w}_f(\alpha) e^{2 \pi i \alpha \cdot x},\] 
where
\[ \hat{w}_f(\alpha) = \int_{\Bb{R}^d} \hat{f}(\xi) \overline{\hat{f}(\xi + \alpha)} \sum_{p \in \Cal{P}_{\alpha}} |\det C_p|^{-1} \overline{\hat{h}_p(\xi)}\hat{g}_p(\xi + \alpha) d \xi.\]
By \cite[Lemma 2.5]{HLW02a} and the proof of Theorem 2.1 in \cite{HLW02a}, $w_f(x)$ is constant for all $f \in \Cal{D}$ if and only if for all $\alpha \in \Lambda \setminus \{0\}$,
\[ \sum_{p \in \Cal{P}_{\alpha}} |\det C_p|^{-1} \overline{\hat{h}_p(\xi)}\hat{g}_p(\xi + \alpha) = 0 a.e. \ \xi.\]

It is well known that if $\Theta$ commutes with $T_{\beta}$ for all $\beta \in \Bb{R}^d$, then it is a Fourier multiplier.  Evaluating $w_f(x)$ at $x = 0$ yields
\[ w_f(0) = \int_{\Bb{R}^d} \hat{f}(\xi) \overline{\hat{f}(\xi)} \sum_{p \in \Cal{P}} |\det C_p|^{-1} \overline{\hat{h}_{p}(\xi)}\hat{g}_{p}(\xi) d \xi = \langle \Theta f, f \rangle.\]
Therefore, since this is valid for all $f \in \Cal{D}$, the symbol of $\Theta$ is $s(\xi)$ as above.
\end{proof}

\begin{corollary} \label{C:Trans1}
Let $\{T_{C_p k} g_p\}$, $\{T_{C_p k} h_p\}$ and $\Theta$ be as in Theorem \ref{T:g_h_alpha}.  We have $\Theta = 0$ if and only if $\Theta \in \{T_{\beta}: \beta \in \Bb{R}^d\}'$ and
\[ s(\xi) = \sum_{p \in \Cal{P}} |\det C_p|^{-1} \overline{\hat{h}_{p}(\xi)} \hat{g}_{p}(\xi) = 0 \ a.e. \xi. \]
Equivalently, $\Theta = 0$ if and only if for each $\alpha \in \Lambda \setminus \{0\}$, equation \ref{E:g_h_alpha} is satisfied and $s(\xi) = 0 $ a.e. $\xi$.
\end{corollary}
\begin{proof}
Clearly, if $\Theta = 0$, then $\Theta \in \{T_{\beta}: \beta \in \Bb{R}^d\}^{\prime}$, whence for all $\alpha \in \Lambda \setminus \{0\}$, equation (\ref{E:g_h_alpha}) is satisfied.  Moreover, $s(\xi) = 0$.  Conversely, if for all $\alpha \in \Lambda \setminus \{0\}$, equation (\ref{E:g_h_alpha}) is satisfied, then $\Theta \in \{T_{\beta}: \beta \in \Bb{R}^d\}^{\prime}$, and if $s(\xi) = 0$ as well, then $\Theta = 0$.
\end{proof}

\subsection{Different Translation Lattices}

\begin{lemma} \label{L:CD1}
Let $\Bb{G} := \{T_{C k} g_p\}$ and $\Bb{H} := \{T_{D k} h_p\}$ be Bessel, and define
\[ \THG:= \sum_{p \in \Cal{P}} \sum_{k \in \Bb{Z}^d} \langle \cdot , T_{C k} g_p \rangle T_{D k} h_p. \]
For all $z \in \Bb{Z}^d$, $\THG T_{Cz} = T_{Dz} \THG$.
\end{lemma}
\begin{proof}
The proof is a simple computation:
\begin{align*}
\THG T_{Cz} &= \sum_{p \in \Cal{P}} \sum_{k \in \Bb{Z}^d} \langle T_{Cz} \cdot , T_{C k} g_p \rangle T_{D k} h_p \\
&= \sum_{p \in \Cal{P}} \sum_{k \in \Bb{Z}^d} \langle \cdot , T_{C (k-z)} g_p \rangle T_{D k} h_p \\
&= \sum_{p \in \Cal{P}} \sum_{k \in \Bb{Z}^d} \langle \cdot , T_{C k} g_p \rangle T_{D (k + z)} h_p \\
&= T_{D z} \THG.
\end{align*}
Since the sequence $\{T_{Dk}g_p\}$ is Bessel, the sum converges in norm, whence the $T_{D z}$ factors out of the sum.
\end{proof}

\begin{proposition} \label{P:CD1}
Let $\Bb{G}$, $\Bb{H}$, and $\THG$ be as in Lemma \ref{L:CD1}.  If $C \neq D$ and $\THG \in \{T_{C z}: z \in \Bb{Z}^d\}^{\prime}$, then $\THG = 0$.
\end{proposition}
\begin{proof}
If $\THG \in \{T_{C z}: z \in \Bb{Z}^d\}^{\prime}$, then we have by Lemma \ref{L:CD1} that $T_{C z} \THG = \THG T_{C z} = T_{D z} \THG$.  Hence, if $\THG \neq 0$, there exists a non-zero function $f \in \ltwod$ such that $T_{C z}f = T_{D z} f$, hence $T_{Cz - Dz} f = f$.  However, it is well known that if $Cz - Dz \neq 0$, $T_{Cz - Dz}$ has purely continuous spectrum and hence no non-zero eigenvectors.  Therefore, $T_{C z} = T_{D z}$ for all $z \in \Bb{Z}^d$, whence $C = D$.
\end{proof}

\begin{corollary}
If $C \neq D$, then $\Bb{G}$ and $\Bb{H}$ cannot be dual frames.
\end{corollary}
\begin{proof}
If $\Bb{G}$ and $\Bb{H}$ are dual frames, then $\THG = I$, but by Proposition \ref{P:CD1}, this is not possible.
\end{proof}

\begin{proposition} \label{P:TCD}
Let $\Bb{G}$ and $\Bb{H}$ be as in Lemma \ref{L:CD1} and satisfy the local integrability condition, and let $\THG$ be as in Lemma \ref{L:CD1}.  We have $\THG = 0$ if and only if
\[ \sum_{p \in \Cal{P}}  \overline{\hat{g}_{p}(C'\xi)} \hat{h}_{p}(D'(\xi + k)) = 0 \ a.e. \xi \]
for all $k \in \Bb{Z}^d$.
\end{proposition}

\begin{proof}
Let $D_{C}$ and $D_{D}$ be the (unitary) dilation operators associated to the matrices $C$ and $D$, respectively.  By the polarization identity, $\THG = 0$ if and only if for every $f \in \Cal{D}$, $\langle D_{C}^{-1} \THG D_{D} f, f \rangle = 0$.  Recall the commutation relation $D_{C} T_{C z} = T_{z} D_{C}$.

\begin{align*}
\langle D_{D} \THG D_{C}^{-1} f, f \rangle &= \langle \THG D_{C}^{-1} f, D_{D}^{-1} f \rangle \\ 
&= \sum_{p \in \Cal{P}} \sum_{k \in \Bb{Z}^d} \langle D_{C}^{-1} f, T_{C k} g_p \rangle \langle T_{D k} h_p , D_{D}^{-1} f \rangle \\
&= \sum_{p \in \Cal{P}} \sum_{k \in \Bb{Z}^d} \langle f, D_{C} T_{C k} g_p \rangle \langle D_{D} T_{D k} h_p , f \rangle \\
&= \sum_{p \in \Cal{P}} \sum_{k \in \Bb{Z}^d} \langle f, T_{k} D_{C} g_p \rangle \langle T_{k} D_{D} h_p , f \rangle.
\end{align*}
We apply Corollary \ref{C:Trans1} to the systems $\{T_{k} D_{C} g_p\}$ and $\{T_{k} D_{D} h_{p}\}$ (note that these collections satisfy the local integrability condition with respect to the integer lattice).  Here, $\Lambda = \Bb{Z}^d$ and for each $\alpha \in \Lambda, \Cal{P}_{\alpha} = \Cal{P}$.  Therefore, $\THG = 0$ if and only if for each $k \in \Bb{Z}^d$,
\[ \sum_{p \in \Cal{P}} \overline{\widehat{D_{C} g_p}(\xi)} \widehat{D_{D} h_p}(\xi + k) = 0 \ a.e. \ \xi.\]
Since $\widehat{D_{C}} = D_{C'}$, it now follows that $\THG = 0$ if and only if 
\[ \sum_{p \in \Cal{P}} \overline{\hat{g}_{p}(C' \xi)} \hat{h}_{p}(D'(\xi + k)) = 0 \ a.e. \ \xi\]
for every $k \in \Bb{Z}^d$.
\end{proof}

For singly generated systems, we recover the characterization developed in \cite{ALTW02a}.

\begin{corollary}
Suppose $\{T_{Ck}g : k \in \Bb{Z}^d\}$ and $\{T_{Dk}h : k \in \Bb{Z}^d\}$ are Bessel.  Then
\[ \sum_{k \in \Bb{Z}^d} \langle \cdot , T_{C k} g \rangle T_{D k} h = 0 \]
if and only if
\[  \sum_{k \in \Bb{Z}^d} |\hat{g}(C'(\xi + k))|^2 \cdot \sum_{k \in \Bb{Z}^d} |\hat{h}(D'(\xi + k))|^2 = 0 \ a.e. \xi. \]
\end{corollary}
\begin{proof}
For singly generated systems, the Bessel condition is equivalent to the local integrability condition \cite{ALTW02a}.  If the Bessel sequences are orthogonal, then for each $k \in \Bb{Z}^d$,
\[ \overline{\hat{g}(C'\xi)} \hat{h}(D'(\xi + k)) = 0 \ a.e. \xi, \]
hence for each $m \in \Bb{Z}^d$
\[ |\hat{g}(C'(\xi + m)|^2 \cdot |\hat{h}(D'(\xi + m + k))|^2 = 0 \ a.e. \xi. \]
Summing over $m$ and $k$ yields
\begin{align*}
0 &= \sum_{m \in \Bb{Z}^d} \sum_{k \in \Bb{Z}^d} |\hat{g}(C'(\xi + m)|^2 |\hat{h}(D'(\xi + m + k))|^2 \\
    &= \sum_{m \in \Bb{Z}^d} |\hat{g}(C'(\xi + m)|^2 \sum_{k \in \Bb{Z}^d}  |\hat{h}(D'(\xi + m + k))|^2 \\
    &= \sum_{m \in \Bb{Z}^d} |\hat{g}(C'(\xi + m)|^2 \sum_{k \in \Bb{Z}^d}  |\hat{h}(D'(\xi + k))|^2
\end{align*}
for almost every $\xi$.  The converse follows by reversing the steps above.
\end{proof}

\section{Affine Systems}

The fundamental work of Ron and Shen \cite{RS97b, RS97c} shows an intimate connection between affine and quasi-affine reproducing systems for integer dilations.  Recent work by Labate, Hernandez and Weiss \cite{HLW02a} shows that for non-integer dilations, the analogous results do not necessarily hold (see also \cite{CCMW02a}.  The results that follow also show how the two systems are related in some cases, and not related in others in terms of orthogonal systems (see Example \ref{E:qadiff}).

For a $d \times d$ invertible matrix $A$, let $D_{A}$ denote the unitary operator
\[ D_{A}: \ltwod \to \ltwod : f(\cdot ) \mapsto \sqrt{| \det A|} f(A \cdot ) \]
and let $\widetilde{D}_{A}$ denote the renormalized operator
\[ \widetilde{D}_{A} : \ltwod \to \ltwod : f(\cdot) \mapsto | \det A| f(A \cdot ).\]

The affine and quasi-affine systems, respectively, are as follows.
\[ \Cal{U}_{A,X}(\Psi) := \{ D_A^{n} T_{Xz} \psi_i: n \in \Bb{Z}; \ z \in \Bb{Z}^d; \ \psi_i \in \Psi\} \]

\[ \Cal{U}^{q}_{A,X}(\Psi) := \{ D_A^{n} T_{Xz} \psi_i: n \geq 0; \ z \in \Bb{Z}^d; \ \psi_i \in \Psi\} \cup \{T_{Xz} \widetilde{D}_A^{n} \psi_i: n < 0; \ z \in \Bb{Z}^d; \ \psi_i \in \Psi\} \]
In case $X = I$, we shall write $\Cal{U}_{A}(\Psi)$ and $\Cal{U}^{q}_{A}(\Psi)$.  We will always assume that $\Psi$ and $\Phi$ are finite collections in $\ltwod$.  We say that $A$ is expanding if all eigenvalues of $A$ have modulus strictly greater than 1.  We say that $A$ is integer valued if all entries of $A$ are integers.

Note the following commutation relations:  if $A$ is a $d \times d$ invertible matrix and $\alpha \in \Bb{R}^d$, then
\begin{equation} \label{E:commute1}
D_{A} T_{\alpha} = T_{A^{-1} \alpha} D_{A} \qquad \text{ and } \qquad T_{\alpha} D_{A} = D_{A} T_{A \alpha}.
\end{equation}
If $B$ is also a $d \times d$ invertible matrix, then
\begin{equation} \label{E:commute2}
D_{A} D_{B} = D_{\tilde{B}} D_{A}
\end{equation}
where $\tilde{B} = A^{-1} B A$.  Note also that $\Cal{U}_{A,X}(\Psi)$ is dilation invariant, i.e. $D_{A}\Cal{U}_{A,X}(\Psi) \subset \Cal{U}_{A,X}(\Psi)$.  Moreover, if the lattice $X \Bb{Z}^d$ is invariant under the matrix $A$, then by the commutation relation (\ref{E:commute1}), $\Cal{U}^{q}_{A,X}(\Psi)$ is shift invariant, i.e. $T_{X m} \Cal{U}^{q}_{A,X}(\Psi) \subset \Cal{U}^{q}_{A,X}(\Psi)$.

In order to apply the results of section 1, we will view
\[ \Cal{U}_{A,X}(\Psi) = \{ T_{A^{-n} X z} D_{A}^{n} \psi_i\} \]
and
\[ \Cal{U}^{q}_{A,X}(\Psi) : \{ T_{A^{-n} Xz} D_{A}^{n} \psi_i: n \geq 0; \ z \in \Bb{Z}^d; \ \psi_i \in \Psi\} \cup \{T_{Xz} \widetilde{D}_A^{n} \psi_i: n < 0; \ z \in \Bb{Z}^d; \ \psi_i \in \Psi\}. \]
In both cases, $\Cal{P} = \Bb{Z} \times \{1, \dots, n\}$.  For $A$ an expanding matrix, if the affine system $\Cal{U}_{A,X}(\Psi)$ is Bessel, then it also satisfies the local integrability condition, and likewise for the quasi-affine system \cite{HLW02a}.

\begin{lemma} \label{L:CofV}
Suppose $\Psi \subset \ltwod$; $D_{X} \Cal{U}_{A, X} (\Psi) = \Cal{U}_{\tilde{A}} (D_{X} \Psi)$ and $D_{X} \Cal{U}^{q}_{A, X} (\Psi) = \Cal{U}^{q}_{\tilde{A}} (D_{X} \Psi)$, where $\tilde{A} = X^{-1}AX$.
\end{lemma}

\begin{proof} 
The proof is a consequence of the commutation relations (\ref{E:commute1}) and (\ref{E:commute2}); see\cite{Bow00c}.
\end{proof}

\begin{lemma} \label{L:DT}
If $A$ is an expansive matrix, then $\{D_{A}\}' \cap \{T_{z} : z \in \Bb{Z}^d \}'$ is the von Neumann algebra of Fourier multipliers whose symbol $s(\xi)$ satisfies $s(A^{*}\xi) = s(\xi) \ a.e. \ \xi$.  In other words, $S \in \{D_{A}\}' \cap \{T_{z} : z \in \Bb{Z}^d \}'$ if and only if $\widehat{Sf} (\xi) = s(\xi) \hat{f}(\xi)$ for $s(\cdot) \in L^{\infty}(\Bb{R}^d)$ and $s(A^{*}\xi) = s(\xi) \ a.e. \ \xi$.
\end{lemma}
\begin{proof}
Suppose that $S \in \{D_{A}\}' \cap \{T_{z} : z \in \Bb{Z}^d \}'$.  Note that $D_{A} T_{z} D_{A}^{-1} = T_{A^{-1}z}$, whence $S$ commutes with every operator of the form $T_{A^{n} z}$.  Since $A$ is expansive, the set $\cup_{n \in \Bb{Z}} A^{n} \Bb{Z}^d$ is dense in $\Bb{R}^d$; whence the operators $\{T_{A^{n}z}: n \in \Bb{Z}; \ \Bb{Z}^d \}$ are dense in $\{T_{\beta}: \beta \in \Bb{R}^d\}$ in the weak operator topology.  Therefore, $S \in \{T_{\beta} : \beta \in \Bb{R}^d \}'$, and hence is a Fourier multiplier.  Moreover, since $S \in \{D_{A}\}'$, the symbol of $S$ must satisfy $s(A^{*}\xi) = s(\xi) \ a.e. \ \xi$ since for all $f \in \ltwod$:
\[ s(\xi) \sqrt{|\det A|}^{-1} \hat{f}(A'\xi) = \widehat{S D_{A} f}(\xi) = \widehat{D_{A} S f}(\xi) = \sqrt{|\det A|}^{-1} s(A' \xi) f(A'\xi).\]
The reverse implication now follows by the above computation.
\end{proof}

\begin{lemma} \label{L:tq}
Suppose that $A$ is an expansive integral matrix, and suppose that $\Cal{U}_{A}(\Psi)$ and $\Cal{U}_{A}(\Phi)$ are Bessel sequences.  The following are equivalent:
\begin{enumerate}
\item $\theta^{*}_{\Phi} \theta_{\Psi} \in \{T_{\beta}: \beta \in \Bb{R}^d\}'$;
\item $\theta^{q*}_{\Phi} \theta^{q}_{\Psi} \in \{T_{\beta}: \beta \in \Bb{R}^d\}'$;
\item $q \in \Bb{Z}^d \setminus A^{*}\Bb{Z}^d$,
\[ \sum_{i=1}^{r} \sum_{j = 0}^{\infty} \overline{\hat{\phi}_i(A^{*j} \xi)} \hat{\psi}_i(A^{*j}(\xi + q)) = 0 \ a.e. \ \xi. \]
\end{enumerate}
Moreover, in any of the three cases, the symbol both $\theta^{*}_{\Phi} \theta_{\Psi}$ and $\theta^{q*}_{\Phi} \theta^{q}_{\Psi}$ is
\[ s(\xi) = \sum_{i=1}^{r} \sum_{j \in \Bb{Z}} \overline{\hat{\phi}_i(A^{*j}\xi)} \hat{\psi}_i(A^{*j}\xi) \ a.e. \ \xi.\]
\end{lemma}
\begin{proof}
We apply Theorem \ref{T:g_h_alpha} to the affine systems $\Cal{U}_{A}(\Psi)$ and $\Cal{U}_{A}(\Phi)$, and to the quasi-affine systems $\Cal{U}^{q}_{A}(\Psi)$ and $\Cal{U}^{q}_{A}(\Phi)$.  For the affine systems, $\Cal{P} = \Bb{Z} \times \{1,\dots,r\}$; for $z \in \Bb{Z}$, $C_{z,i} = A^{-z}$; $g_{z,i} = D_{A}^z \psi_i$, $h_{z,i} = D_{A}^z \phi_i$; $\Lambda = \cup_{n \in \Bb{Z}} A^{*n} \Bb{Z}^d$, and if $\alpha \in \Lambda$, then $\alpha = A^{*s} q$ for some $s \in \Bb{Z}$ and some $q \in \Bb{Z}^d \setminus A^{*} \Bb{Z}^d$.  (For the remainder of the proof, we will suppress the index $i$).  We have $\Cal{P}_{\alpha} = \{n: A^{* -n} A^{*s}q \in \Bb{Z}^d\} = \{n: s \geq n\}$.  Therefore,
\begin{align*}
\sum_{p \in \Cal{P}_{\alpha}} |\det C_{p}|^{-1} \overline{\hat{g}_{p}(\xi)}\hat{h}_{p}(\xi + \alpha) &= \sum_{i=1}^{r} \sum_{n = -\infty}^{s} |\det A^{-n}|^{-1} \overline{D_{A}^{-n}\hat{\psi}_i(\xi)} D_{A}^{-n}\hat{\phi}_i(\xi + A^{*s}q) \\
&= \sum_{i=1}^{r} \sum_{n = -\infty}^{s} \overline{\hat{\psi}_i(A^{*-n}\xi)}\hat{\phi}_i(A^{*-n}(\xi + A^{*s}q)) \\
&= \sum_{i=1}^{r} \sum_{n = -\infty}^{s} \overline{\hat{\psi}_i(A^{*-n+s} A^{*-s}\xi)}\hat{\phi}_i(A^{*-n+s}(A^{*-s}\xi + q)) \\
&= \sum_{i=1}^{r} \sum_{n = 0}^{\infty} \overline{\hat{\psi}_i(A^{*n} A^{*-s}\xi)}\hat{\phi}_i(A^{*n}(A^{*-s}\xi + q)).
\end{align*} 

Likewise, for the quasi-affine system, $\Cal{P} = \Bb{Z} \times \{1,\dots,r\}$.  However, for $z > 0$, $C_{z,i} = A^{-z}$ and for $z \leq 0$, $C_{z,i} = I$.  For $z > 0$, $g_{z,i} = D_{A}^z \psi_i$ and $h_{z,i} = D_{A}^z \phi_i$ and for $z \leq 0$, $g_{z,i} = \widetilde{D}_{A}^z \psi_i$ and $h_{z,i} = \widetilde{D}_{A}^{n} \phi_i$.   Here $\Lambda = \Bb{Z}^d$, and if $\alpha \in \Lambda$, then $\alpha = A^{*s} q$ for some $s \in \Bb{Z}$ and some $q \in \Bb{Z}^d \setminus A^{*} \Bb{Z}^d$.  (Again we will suppress the index $i$).  We have $\Cal{P}_{\alpha} = \{n > 0: A^{* -n} A^{*s}q \in \Bb{Z}^d\} \cup \{n: n \leq 0\} = \{n: s \geq n\}$.  Therefore,
\begin{align*}
\sum_{p \in \Cal{P}_{\alpha}} |\det C_{p}|^{-1} \overline{\hat{g}_{p}(\xi)}\hat{h}_{p}(\xi + \alpha) &= \sum_{i=1}^{r} \sum_{n = 1}^{s} |\det A^{-n}|^{-1} \overline{D_{A}^{-n}\hat{\psi}_i(\xi)} D_{A}^{-n}\hat{\phi}_i(\xi + A^{*s}q) \\
&  \qquad \qquad + \sum_{i=1}^{r} \sum_{n = -\infty}^{0} |\det I|^{-1} \overline{\widetilde{D}_{A}^{-n}\hat{\psi}_i(\xi)} \widetilde{D}_{A}^{-n}\hat{\phi}_i(\xi + A^{*s}q) \\
&= \sum_{i=1}^{r} \sum_{n = -\infty}^{s} \overline{\hat{\psi}_i(A^{*-n}\xi)}\hat{\phi}_i(A^{*-n}(\xi + A^{*s}q)) \\
&= \sum_{i=1}^{r} \sum_{n = 0}^{\infty} \overline{\hat{\psi}_i(A^{*n} A^{*-s}\xi)}\hat{\phi}_i(A^{*n}(A^{*-s}\xi + q)).
\end{align*}

The lemma now follows by Theorem \ref{T:g_h_alpha}.
\end{proof}

\begin{theorem} \label{T:AQE}
Suppose $\Cal{U}_{A}(\Psi)$ and $\Cal{U}_{A}(\Phi)$ are Bessel sequences, where $A$ is an expansive integral matrix.  The following are equivalent:
\begin{enumerate}
\item $\theta^{*}_{\Phi} \theta_{\Psi} \in \{T_z: z \in \Bb{Z}^d\}'$;
\item $\theta^{q*}_{\Phi} \theta^{q}_{\Psi} \in \{D_{A}\}'$;
\item $\theta^{*}_{\Phi} \theta_{\Psi} \in \{T_\beta: \beta \in \Bb{R}^d\}'$
\item $\theta^{q*}_{\Phi} \theta^{q}_{\Psi} \in \{T_\beta: \beta \in \Bb{R}^d\}'$
\item $\theta^{*}_{\Phi} \theta_{\Psi} = \theta^{q*}_{\Phi} \theta^{q}_{\Psi}$;
\item for $q \in \Bb{Z}^d \setminus A^{*} \Bb{Z}^d$, $\displaystyle{\sum_{i=1}^{r} \sum_{j = 0}^{\infty} \overline{\hat{\phi}_i(A^{*j} \xi)} \hat{\psi}_i(A^{*j}(\xi + q)) = 0 \ a.e. \ \xi}$;
\item $\theta^{*}_{\Phi} \theta_{\Psi}$ is a Fourier multiplier, i.e. $\widehat{\theta^{*}_{\Phi} \theta_{\Psi} f} (\xi) = s(\xi) \hat{f}(\xi)$, whose symbol is
\[ s(\xi) = \sum_{i=1}^{r} \sum_{j \in \Bb{Z}} \overline{\hat{\phi}_i(A^{*j}\xi)} \hat{\psi}_i(A^{*j}\xi) \ a.e. \ \xi;\]
\item $\theta^{q*}_{\Phi} \theta^{q}_{\Psi}$ is a Fourier multiplier, with the same symbol $s(\xi)$.
\end{enumerate}
\end{theorem}
\begin{proof}
The implications $1 \Rightarrow 3$ and $2 \Rightarrow 4$ follow from Lemma \ref{L:DT}.  The symbol $s(\xi)$ above satisfies $s(A^{*}\xi) = s(\xi) \ a.e. \ \xi$, hence the implications $7 \Rightarrow 1$ and $8 \Rightarrow 2$ also follow from Lemma \ref{L:DT}.  Lemma \ref{L:tq} yields $3 \Rightarrow 6$, $4 \Rightarrow 6$, $6 \Rightarrow 7$, $6 \Rightarrow 8$, and $5 \Leftrightarrow 6$.

Thus we have demonstrated
\[ 7 \Rightarrow 1 \Rightarrow 3 \Rightarrow 6 \Rightarrow 8 \Rightarrow 2 \Rightarrow 4 \Rightarrow 6 \Rightarrow 7 \text{ and } 5 \Leftrightarrow 6. \]
\end{proof}

\begin{remark} If $\Psi = \Phi$ in the preceeding theorem, the conditions there are equivalent to the condition that the canonical dual of $\Cal{U}_{A}^{q}(\Psi)$ also has the quasi-affine structure \cite{BW03a}.
\end{remark}

\begin{corollary} \label{C:OWF}
Suppose $A$ is an expansive integral matrix and suppose $\Cal{U}_{A}(\Psi)$ and $\Cal{U}_{A}(\Phi)$ are Bessel sequences.  Then they are orthogonal if and only if for all $q \in \Bb{Z}^d \setminus A^{*}\Bb{Z}^d$,
\[\sum_{i=1}^{r} \sum_{j \geq 0} \overline{\hat{\psi}_i(A^{*j} \xi)} \hat{\phi}_i(A^{*j}(\xi + q)) = 0 \ a.e. \ \xi;\]
and 
\[ \sum_{i=1}^{r} \sum_{j \in \Bb{Z}} \overline{\hat{\psi}_i(A^{*j} \xi)} \hat{\phi}_i(A^{*j}\xi) = 0 \ a.e. \ \xi. \]
Moreover, the Bessel sequences $\Cal{U}^{q}_{A}(\Psi)$ and $\Cal{U}^{q}_{A}(\Phi)$ are orthogonal if and only if the same two equations hold.
\end{corollary} 
\begin{proof}
This follows immediately from Theorem \ref{T:AQE}.
\end{proof}

We now consider the case when the two affine systems have different dilation matrices and/or different translation lattices.

\begin{lemma} \label{L:AAB}
Suppose $\Cal{U}_{A}(\Psi)$ and $\Cal{U}_{B}(\Phi)$ are Bessel.  If $\theta_{\Psi}^{*} \theta_{\Phi} \in \{T_{\beta}: \beta \in \Bb{R}^d\}'$, then either: 1) $A=B$ or 2) $\theta_{\psi}^{*} \theta_{\phi} = 0$.
\end{lemma}
\begin{proof}
See Proposition \ref{P:CD1}.
\end{proof}

\begin{proposition}
If the frames $\Cal{U}_{A,X}(\Psi)$ and $\Cal{U}_{B,Y}(\Phi)$ are dual, then $A=B$.
\end{proposition}
\begin{proof}
This follows directly from Lemma \ref{L:AAB}.
\end{proof}

\begin{lemma} \label{L:thetas}
Let $A$ and $B$ be any expansive matrices and suppose that $\Cal{U}_{A}(\Psi)$ and $\Cal{U}_{B}(\Phi)$ are Bessel sequences.  The following are equivalent:
\begin{enumerate}
\item $\Theta_{A,B}^{+} := \sum_{i = 1}^{r} \sum_{n > 0} \sum_{z \in \Bb{Z}^d} \langle \cdot , D_{A}^n T_z \psi_i \rangle D_{B}^n T_z \phi_i = 0$;
\item $\Theta_{A,B}^{-} := \sum_{i = 1}^{r} \sum_{n < 0} \sum_{z \in \Bb{Z}^d} \langle \cdot , D_{A}^n T_z \psi_i \rangle D_{B}^n T_z \phi_i = 0$;
\item $\Theta_{A,B}^{0} := \sum_{i = 1}^{r} \sum_{z \in \Bb{Z}^d} \langle \cdot , T_z \psi_i \rangle T_z \phi_i = 0$;
\item for all $k \in \Bb{Z}^d$, $\sum_{i = 1}^{r} \hat{\psi}_i(\xi) \overline{\hat{\phi}_i(\xi + k)} = 0 \ a.e. \ \xi$.
\end{enumerate}
\end{lemma}
\begin{proof}
The equivalence of items 3 and 4 follow from Proposition \ref{P:TCD}.  Note that $\Theta_{A,B}^{+} = \sum_{n > 0} D_{B}^n \Theta_{A,B}^{0} D_{A}^{-n}$, and similarly for $\Theta_{A,B}^{-}$, hence item 3 implies items 1 and 2.  Consider the following computation:
\begin{align} \label{E:thetas}
D_{B}^{-1} \Theta_{A,B}^{+} D_{A} &= D_{B}^{-1} \sum_{i = 1}^{r} \sum_{n > 0} \sum_{z \in \Bb{Z}^d} \langle D_{A} \cdot , D_{A}^n T_z \psi_i \rangle D_{B}^n T_z \phi_i \\
&= \sum_{i = 1}^{r} \sum_{n > 0} \sum_{z \in \Bb{Z}^d} \langle \cdot , D_{A}^{n-1} T_z \psi_i \rangle D_{B}^{n-1} T_z \phi_i \notag \\
&= \Theta_{A,B}^{+} + \Theta_{A,B}^{0}. \notag
\end{align}
Therefore, 1 implies 3.  An analogous computation shows 2 implies 3.
\end{proof}

\begin{lemma} \label{L:AisB}
Suppose $A$ and $B$ are expansive matrices, with $A$ integer valued, and suppose that $\Cal{U}_{A}(\Psi)$ and $\Cal{U}_{B}(\Phi)$ are Bessel sequences.  Let $\Theta_{A,B}^{+}$ be as in Lemma \ref{L:thetas}.  If $\Theta_{A,B}^{+} \in \{T_z: z \in \Bb{Z}^d\}'$, then either $A=B$ or $\Theta_{A,B}^{+} = 0$.
\end{lemma}
\begin{proof}
By the computation in the proof of Lemma \ref{L:thetas}, we have 
\[ D_{B}^{-1} \Theta_{A,B}^{+} D_{A} = \Theta_{A,B}^{+} + \Theta_{A,B}^{0}.\]
Since $\Theta_{A,B}^{0} \in \{T_z: z \in \Bb{Z}^d\}'$, if $\Theta_{A,B}^{+} \in \{T_z: z \in \Bb{Z}^d\}'$, then $D_{B}^{-1} \Theta_{A,B}^{+} D_{A} \in \{T_z: z \in \Bb{Z}^d\}'$ as well.  Therefore for all $z \in \Bb{Z}^d$, by the commutation relation for translations and dilations,
\[ D_{B}^{-1} \Theta_{A,B}^{+} D_{A} T_{A z} = D_{B}^{-1} \Theta_{A,B}^{+} T_{z} D_{A} = D_{B}^{-1} T_{z} \Theta_{A,B}^{+} D_{A} = T_{Bz} D_{B}^{-1} \Theta_{A,B}^{+} D_{A}. \]
Therefore, 
\[ T_{Bz} D_{B}^{-1} \Theta_{A,B}^{+} D_{A} = D_{B}^{-1} \Theta_{A,B}^{+} D_{A} T_{Az} = T_{Az} D_{B}^{-1} \Theta_{A,B}^{+} D_{A}.\]
Hence, if $\Theta_{A,B}^{+} \neq 0$, then there exists a function $f \in \ltwod$ such that $T_{Bz}f = T_{Az}f$ for all $z \in \Bb{Z}^d$.  It follows that $Bz = Az$ and hence $A = B$.
\end{proof}

We end this subsection with the following result, which is not a complete characterization but the best possible result with the present techniques.
\begin{theorem} \label{T:SXY}
A sufficient condition for the Bessel sequences $\Cal{U}_{A,X}(\Psi)$ and $\Cal{U}_{B,Y}(\Phi)$ to be orthogonal is
\[ \sum_{i =1}^{r} \overline{\hat{\psi}_i(X' \xi)} \hat{\phi}_i(Y'(\xi + k)) = 0 \ a.e.,\]
for all $k \in \Bb{Z}^d$.
\end{theorem}
\begin{proof}
If $\sum_{i =1}^{r} \overline{\hat{\psi}_i(X' \xi)} \hat{\phi}_i(Y'(\xi + k)) = 0 \ a.e.$ for all $k \in \Bb{Z}^d$, then by Proposition \ref{P:TCD}, $\{T_{Xz} \psi_i\}$ and $\{T_{Yz} \phi_i\}$ are orthogonal.  It follows then by Lemma \ref{L:thetas} that the affine sets are orthogonal.
\end{proof}

\subsection{Quasi-Affine Systems}

\begin{theorem} \label{T:quasiAB}
Let $A$ be an expansive integral matrix, and let $B$ be any expansive matrix such that $A \neq B$.  Suppose the quasi-affine systems $\Cal{U}^{q}_{A}(\Psi)$ and $\Cal{U}^{q}_{B}(\Phi)$ are Bessel; they are orthogonal if and only if
\begin{enumerate}
\item $\sum_{i=1}^{r} \overline{\hat{\psi}_{i}(\xi)} \hat{\phi}_{i}(\xi + k) = 0 \ a.e. \ \xi$ for every $k \in \Bb{Z}^d$;
\item $\sum_{i=1}^{r} \sum_{j > 0} \overline{\hat{\psi}_{i}(A^{*j} \xi)} \hat{\phi}_{i}(B^{*j}(\xi + k)) = 0 \ a.e. \ \xi$ for all $k \in \Bb{Z}^d$.
\end{enumerate}
\end{theorem}

\begin{proof}
Write the operator $\theta_{\Phi}^{q *} \theta_{\Psi}^{q}$ as the sum $M + N$, where 
\[ M: = \sum_{i=1}^{r} \sum_{n < 0} \sum_{z \in \Bb{Z}^d} \langle \cdot , T_z \tilde{D}_{A}^{n} \psi_i \rangle T_z \tilde{D}_{B}^{n} \phi_i \ \text{ and } \ 
N:= \sum_{i=1}^{r} \sum_{n \geq 0} \sum_{z \in \Bb{Z}^d} \langle \cdot , D_{A}^{n} T_z \psi_i \rangle D_{B}^{n} T_z \phi_i. \]
By definition, $M \in \{ T_z : z \in \Bb{Z}^d\}'$, thus if $\theta_{\Phi}^{q *} \theta_{\Psi}^{q} = M + N = 0$ and $A \neq B$, then by Lemma \ref{L:AisB}, $N = 0$.  Therefore, $\theta_{\Phi}^{q *} \theta_{\Psi}^{q} = 0$ if and only if $M = N = 0$.

By Lemma \ref{L:thetas}, $N = 0$ if and only if item 1.  By Theorem \ref{T:g_h_alpha} and Corollary \ref{C:Trans1}, $M = 0$ if and only if for each $k \in \Bb{Z}^d$,
\[ \sum_{i = 1}^{r} \sum_{j < 0} \overline{ \widehat{\tilde{D}^{j}_{A} \psi_i}(\xi)} \widehat{\tilde{D}^{j}_{B} \phi_i}(\xi + k) =\sum_{i = 1}^{r} \sum_{j > 0} \overline{ \hat{\psi}_i(A^{*j} \xi)} \hat{\phi}_i(B^{*j} (\xi + k))= 0 \]
for almost every $\xi$.
\end{proof}

\begin{corollary}
If the quasi-affine frames $\Cal{U}^{q}_{A}(\Psi)$ and $\Cal{U}^{q}_{B}(\Phi)$ are dual, then $A=B$.
\end{corollary}
\begin{proof}
Let $M,N$ be as in the proof of Theorem \ref{T:quasiAB}.  If $\theta_{\Phi}^{q *} \theta_{\Psi}^{q} = I$, then $N \in \{ T_z : z \in \Bb{Z}^d\}'$, whence by Lemma \ref{L:AisB}, $A = B$.
\end{proof}

The following corollary is nearly a complete characterization of when quasi-affine systems are orthogonal.

\begin{corollary} \label{C:QAX}
Let $A$ and $B$ be any expansive matrices and $X$ and $Y$ be invertible matrices such that $\tilde{A} := X^{-1}A X$ is an integer matrix and $X^{-1}A X \neq Y^{-1} B Y =: \tilde{Y}$.  Suppose the quasi-affine systems $\Cal{U}^{q}_{A,X}(\Psi)$ and $\Cal{U}^{q}_{B,Y}(\Phi)$ are Bessel; they are orthogonal if and only if
\begin{enumerate}
\item $\sum_{i=1}^{r} \overline{\hat{\psi}_{i}(X' \xi)} \hat{\phi}_{i}(Y'(\xi + k)) = 0 \ a.e. \ \xi$ for every $k \in \Bb{Z}^d$;
\item $\sum_{i=1}^{r} \sum_{j > 0} \overline{\hat{\psi}_{i}(A^{*j} X' \xi)} \hat{\phi}_{i}(B^{*j}Y'(\xi + k)) = 0 \ a.e. \ \xi$ for all $k \in \Bb{Z}^d$.
\end{enumerate}
\end{corollary}

\begin{proof}
The quasi-affine systems $\Cal{U}^{q}_{A,X}(\Psi)$ and $\Cal{U}^{q}_{B,Y}(\Phi)$ are orthogonal if and only if $D_{X} \Cal{U}^{q}_{A,X}(\Psi)$ and $D_{Y} \Cal{U}^{q}_{B,Y}(\Phi)$ are orthogonal.  By Lemma \ref{L:CofV},
\[ D_{X} \Cal{U}^{q}_{A,X}(\Psi) = \Cal{U}^{q}_{\tilde{A}}(D_{X} \Psi) \text{ and } D_{Y} \Cal{U}^{q}_{B,Y}(\Phi) = \Cal{U}^{q}_{\tilde{B}}(D_{Y} \Phi).\]
By Theorem \ref{T:quasiAB}, $\Cal{U}^{q}_{\tilde{A}}(D_{X} \Psi)$ and $\Cal{U}^{q}_{\tilde{B}}(D_{Y} \Phi)$ are orthogonal if and only if for every $k \in \Bb{Z}^d$ and almost every $\xi$,
\[ 0 = \sum_{i = 1}^{r} \overline{ \widehat{D_{X} \psi_i}(\xi)} \widehat{D_{Y} \phi_i}(\xi + k) = \sum_{i=1}^{r} \sqrt{ | \det X Y |}^{-1} \overline{ \hat{\psi}_i(X' \xi)} \hat{\phi}_i(Y' (\xi + k));\]
and
\[ 0 = \sum_{i=1}^{r} \sum_{j < 0} \overline{ \widehat{ \tilde{D}^{j}_{\tilde{A}} D_{X} \psi_i}(\xi)} \widehat{ \tilde{D}^{j}_{\tilde{B}} D_{Y} \phi_i}(\xi + k) = \sum_{i=1}^{r} \sum_{j > 0} \sqrt{| \det XY |}^{-1} \overline{ \hat{\psi}_i(X' \tilde{A}^{*j} \xi)} \hat{\phi}_i(Y' \tilde{B}^{*j} (\xi + k)).\]
However, $X' \tilde{A}^{*j} = A^{*j} X'$ and $Y' \tilde{B}^{*j} = B^{*j} Y'$, so we have
\[ 0 = \sum_{i=1}^{r} \sum_{j > 0} \overline{ \hat{\psi}_i(A^{*j} X' \xi)} \hat{\phi}_i(B^{*j} Y'(\xi + k)).\]
\end{proof}

\begin{corollary}
Let $A$,$B$,$X$,$Y$ be as in Corollary \ref{C:QAX}.  If the quasi-affine Bessel systems $\Cal{U}^{q}_{A,X}(\Psi)$ and $\Cal{U}^{q}_{B,Y}(\Phi)$ are orthogonal, then the affine Bessel systems $\Cal{U}_{A,X}(\Psi)$ and $\Cal{U}_{B,Y}(\Phi)$ are also orthogonal.
\end{corollary}
\begin{proof}
By item 1. in Corollary \ref{C:QAX} and Theorem \ref{T:SXY}, the affine systems are orthogonal.
\end{proof}

\begin{example} \label{E:qadiff}
The following example demonstrates that when the dilations are different, it is possible for the affine systems to be orthogonal while the quasi-affine systems are not.  Let $\psi$ be a Frazier-Jawerth frame wavelet, i.e.~such that $\hat{\psi}$ is symmetric, non-negative, supported on $[-1/32,-1/128] \cup [1/128, 1/32]$ and such that $\sum_{j} \hat{\psi}(2^j \xi) \equiv 1$ (see \cite{FJ90a}).  Now, let $\phi$ be a Frazier-Jawerth frame wavelet for dilation by 3 such that $\hat{\phi}$ is symmetric, non-negative, supported on $[-1/3,-1/27] \cup [1/27,1/3]$ and such that $\sum_{j} \hat{\phi}(3^j \xi) \equiv 1$.  Therefore, $\Cal{U}_{2}(\psi)$, $\Cal{U}_{2}^{q}(\psi)$, $\Cal{U}_{3}(\phi)$, and $\Cal{U}_{3}^{q}(\phi)$ are all Parseval frames for $\ltwo$.

Clearly for all $k \in \Bb{Z}$ we have $\overline{\hat{\psi}(\xi)}\hat{\phi}(\xi + k)$, whence by Theorem \ref{T:SXY}, $\Cal{U}_{2}(\psi)$ and $\Cal{U}_{3}(\phi)$ are orthogonal.  However, since both $\hat{\psi}$ and $\hat{\phi}$ are non-negative, 
\[ \sum_{j > 0} \overline{\hat{\psi}(2^{*j} \xi)} \hat{\phi}(3^{*j}(\xi + k)) \neq 0 \] 
on a set of positive measure, whence by Theorem \ref{T:quasiAB}, $\Cal{U}_{2}^{q}(\psi)$ and $\Cal{U}_{3}^{q}(\phi)$ are not orthogonal.
\end{example}

\subsection{Super-Wavelets}

Super-wavelets were introduced in \cite{HL00a}.  The idea of super frames was also studied in \cite{B99a} in the case of Weyl-Heisenberg frames.  Consider the Hilbert space $\ltwod \oplus \ltwod \oplus \cdots \oplus \ltwod$, the direct sum of $\ltwod$ $r$ times.  Denote this space by $\ltwod^r$.  Define the translation and dilation operators $\overline{T}_z$ and $\overline{D}_{A}$ on $\ltwod^r$ by $\overline{T}_z = T_z \oplus T_z \oplus \cdots \oplus T_z$ and $\overline{D}_{A} = D_{A} \oplus D_{A} \oplus \cdots \oplus D_{A}$.  A (orthonormal) superwavelet is a vector $\Psi = \psi_{1} \oplus \psi_{2} \oplus \cdots \oplus \psi_{r} \in \ltwod^r$ such that
\[ \overline{\Cal{U}}_{A}(\Psi) := \{\overline{D}_{A}^{k} \overline{T}_z \Psi: k \in \Bb{Z} \ z \in \Bb{Z}^d \} \]
is an orthonormal basis of $\ltwod^r$.  A complete characterization of orthonormal superwavelets is obtained in \cite{HL00a}.  
\begin{definition}
A Parseval superwavelet is a vector of the form $\Psi = \psi_{1} \oplus \psi_{2} \oplus \cdots \oplus \psi_{r} \in \ltwod^r$ such that
\[ \overline{\Cal{U}}_{A}(\Psi) := \{\overline{D}_{A}^{k} \overline{T}_z \Psi: k \in \Bb{Z} \ z \in \Bb{Z}^d \} \]
is a Parseval frame of $\ltwod^r$.
\end{definition}
We present below a characterization of Parseval superwavelets.

\begin{theorem}
Suppose $\Cal{U}(\psi_i)$ are Bessel sequences for $i = 1,\dots, r$, and $A$ is an expansive integral matrix.  The following are equivalent:
\begin{enumerate}
\item $\overline{\Cal{U}}_{A}(\psi_1 \oplus \dots \oplus \psi_r)$ is a Parseval frame for $L^2(\Bb{R}^d)^n$;
\item $\overline{\Cal{U}}^{q}_{A}(\psi_1 \oplus \dots \oplus \psi_r)$ is a Parseval frame for $L^2(\Bb{R}^d)^n$;
\item the following equations are satisfied
\begin{enumerate}
\item $\sum_{n \in \Bb{Z}} \hat{\psi}_i(A^{*n} \xi) \overline{\hat{\psi}_j(A^{*n} \xi)} = \delta_{i,j} \ a.e \ \xi$ for $i,j=1,\dots,r$, and
\item $\sum_{n = 0}^{\infty} \hat{\psi}_i(A^{*n} \xi) \overline{\hat{\psi}_j(A^{*n} (\xi + k))} = 0 \ a.e \ \xi$ for $k \in \Bb{Z}^d \setminus A^{*}\Bb{Z}^d$ and $i,j=1,\dots,r$.
\end{enumerate}
\end{enumerate}
\end{theorem}

\begin{proof}
Without loss of generality, assume $r = 2$.  Suppose that $\overline{\Cal{U}}_{A}(\psi_1 \oplus \psi_2)$ is a Parseval frame for $\ltwod^2$.  Let $P$ be the orthogonal projection onto the first coordinate of $\ltwod^2$.  By definition of $\overline{D}_{A}$ and $\overline{T}_z$, both are in $\{P\}'$.  A straight forward computation shows that $\Cal{U}_{A}(\psi_1)$ is a Parseval frame for $\ltwod$, since it is the image of $\overline{\Cal{U}}_{A}(\psi_1 \oplus \psi_2)$ under the projection $P$ (see \cite{A95a}).  Note that since $\overline{\Cal{U}}_{A}(\psi_1 \oplus \psi_2)$ is Parseval,
\[ \Theta := \sum_{k \in \Bb{Z}} \sum_{z \in \Bb{Z}^d} \langle \cdot , \overline{D}_{A} \overline{T}_z \psi_1 \oplus \psi_2 \rangle \overline{D}_{A} \overline{T}_z \psi_1 \oplus \psi_2 = I.\]
Therefore, since $\Cal{U}_{A}(\psi_2)$ is the image of $\overline{\Cal{U}}_{A}(\psi_1 \oplus \psi_2)$ under the projection $P^{\perp}$ and $P$ commutes with $\Theta$, by Lemma \ref{L:PXY}, $\Cal{U}_{A}(\psi_1)$ and $\Cal{U}_{A}(\psi_2)$ are orthogonal.  Combining the characterization theorem for Parseval wavelet frames \cite{HW} with Corollary \ref{C:OWF}, we see that item 1.~implies item 3.

Conversely, if $\Cal{U}_{A}(\psi_1)$ and $\Cal{U}_{A}(\psi_2)$ are both Parseval and are orthogonal, then $\overline{\Cal{U}}_{A}(\psi_1 \oplus \psi_2)$ is also Parseval \cite[Theorem 2.9]{HL00a}, thus item 3.~implies item 1.

The equivalence of items 2.~and 3.~are completely analogous.
\end{proof}

\begin{corollary}
Suppose $A_i$ are (different) expansive integral matrices, and suppose that $\Cal{U}_{A_i}(\psi_i)$ are Parseval frames.  Then $\overline{\Cal{U}}^{q}_{A}(\psi_1 \oplus \dots \oplus \psi_r)$ is a Parseval frame for $L^2(\Bb{R}^d)^r$ if and only if for $i \neq j$ and $k \in \Bb{Z}^d$,
\begin{enumerate}
\item $\overline{\hat{\psi}_i(\xi)} \hat{\psi}_{j}(\xi + k) = 0 \ a.e. \ \xi$;
\item $\sum_{n > 0} \overline{ \hat{\psi}_i(A_{i}^{*n} \xi)} \hat{\psi}_j(A_j^{*n}(\xi + k)) = 0 \ a.e. \ \xi$.
\end{enumerate}
Moreover, if $\overline{\Cal{U}}^{q}_{A}(\psi_1 \oplus \dots \oplus \psi_r)$ is a Parseval frame for $L^2(\Bb{R}^d)^r$, then $\overline{\Cal{U}}_{A}(\psi_1 \oplus \dots \oplus \psi_r)$ is a Parseval frame for $L^2(\Bb{R}^d)^r$.
\end{corollary}

Here we let $\overline{D}_{A} := D_{A_1} \oplus D_{A_2} \oplus \dots \oplus D_{A_r}$ and $\overline{\Cal{U}}^{q}_{A}(\psi_1 \oplus \dots \oplus \psi_r)$ is as before with this dilation operator.

\begin{remark}
Example \ref{E:qadiff} shows that it is possible for $\Cal{U}_{A}(\psi_1 \oplus \dots \oplus \psi_r)$ to be a Parseval superwavelet while $\Cal{U}^{q}_{A}(\psi_1 \oplus \dots \oplus \psi_r)$ is not.
\end{remark}

\section{Dual Frames for Translation Invariant Subspaces} \label{S:DFTIS}

A subspace $M \subset \ltwod$ is translation invariant if for every $\beta \in \Bb{R}^d$, $T_{\beta} M \subset M$.  This is equivalent to the existence of some measurable set $E \subset \Bb{R}^d$ such that
\[ M = \{ f \in \ltwod: supp{\hat{f}} \subset E \}. \]
If $M$ is translation invariant, denote it by $V_{E}$.

\begin{proposition}
Let $\{T_{C_p k} g_p\}$ and $\{T_{C_p k} h_p\}$ be Bessel and satisfy the local integrability condition, and let $E \subset \ltwod$ be measurable.  If equation (\ref{E:g_h_alpha}) is satisfied for every $\alpha \in \Lambda \setminus \{0\}$ and
\[ s(\xi) = \sum_{p \in \Cal{P}} | \det C_{p}|^{-1}\overline{\hat{h}_{p}(\xi)} \hat{g}_{p}(\xi) = 1 \ a.e. \ \xi \in E,\]
then $\{T_{C_p k} g_p\}$ and $\{T_{C_p k} h_p\}$ are $V_{E}$-subspace dual frames.
\end{proposition}
\begin{proof}
By Theorem \ref{T:g_h_alpha}, 
\[ \Theta = \sum_{p \in \Cal{P}} \sum_{k \in \Bb{Z}^d} \langle \cdot , T_{C_p k} g_p \rangle T_{C_p k} h_p \]
is a Fourier multiplier whose symbol is identically 1 on $E$.  It follows that for all $v \in V_E$,
\[ v = \sum_{p \in \Cal{P}} \sum_{k \in \Bb{Z}^d} \langle v , T_{C_p k} g_p \rangle T_{C_p k} h_p. \]
\end{proof}

\begin{proposition}
Suppose $\{T_{C k} g_p\}$ and $\{T_{D k} h_p\}$ are Bessel sequences and let $E$ be measurable.  If $C \neq D$, then $\{T_{D k} h_p\}$ cannot be a $V_{E}$-subspace dual to $\{T_{C k} g_p\}$.
\end{proposition}
\begin{proof}
Let $P$ be the orthogonal projection onto $V_{E}$.  If $\{T_{D k} h_p\}$ is a $V_{E}$-subspace dual to $\{T_{C k} g_p\}$, then by Lemma \ref{L:MSDF}
\[ P = \sum_{p \in \Cal{P}} \sum_{k \in \Bb{Z}^d} \langle \cdot , T_{C_p k} P g_p \rangle T_{D_p k} P h_p \in \{T_{\beta}: \beta \in \Bb{R}^d \}' ,\]
however, by the computation in Proposition \ref{P:CD1}, this is not possible if $C \neq D$.
\end{proof}

\begin{theorem}
Let $\{T_{C_p k} g_p\}$ and $\{T_{C_p k} h_p\}$ be Bessel and satisfy the local integrability condition, and let $E \subset \ltwod$ be measurable.  Then $\{T_{C_p k} h_p\}$ is a $V_E$-subspace dual frame to $\{T_{C_p k} g_p\}$ if and only if
\begin{enumerate}
\item $\displaystyle{\sum_{p \in \Cal{P}} | \det C_{p} |^{-1} \overline{\hat{h}_{p}(\xi)} \hat{g}_{p}(\xi) = 1 \ a.e. \ \xi \in E;}$
\item for all $\alpha \in \Lambda \setminus \{0\}$,
$\displaystyle{\sum_{p \in \Cal{P}_{\alpha}} | \det C_{p} |^{-1} \overline{\hat{h}_{p}(\xi - \alpha)} \hat{g}_{p}(\xi) = 0 \ a.e. \ \xi \in E.}$
\end{enumerate}
\end{theorem}
\begin{proof}
We apply Lemma \ref{L:MSDF} to $\{T_{C_p k} g_p\}$ and $\{T_{C_p k} h_p\}$.  Let $P$ be the projection onto $V_E$; note that $\hat{P} = M_{\chi_{E}}$, i.e.~multiplication by the characteristic function of $E$.  Note also that $P T_{C_p k} = T_{C_p k}P$.

By Lemma \ref{L:MSDF}, we must have for all $v \in V_{E}$,
\[ v = \sum_{p \in \Cal{P}} \sum_{k \in \Bb{Z}^d} \langle v, T_{C_p k} P g_p \rangle T_{C_p k} P h_p, \]
which is equivalent to
\[ \sum_{p \in \Cal{P}} \sum_{k \in \Bb{Z}^d} \langle \cdot , T_{C_p k} P g_p \rangle T_{C_p k} P h_p = P \]
since for all $w \in V_{E}^{\perp}$,
\[ \sum_{p \in \Cal{P}} \sum_{k \in \Bb{Z}^d} \langle w, T_{C_p k} P g_p \rangle T_{C_p k} P h_p = 0.\]
Thus, by Theorem \ref{T:g_h_alpha}, for every $\alpha \in \Lambda$, we must have for almost every $\xi$:
\begin{align} \label{E:local1}
\delta_{\alpha} \chi_{E}(\xi) &= \sum_{p \in \Cal{P}_{\alpha}} | \det C_{p} |^{-1} \overline{\widehat{P h}_{p}(\xi)} \widehat{P g}_{p}(\xi + \alpha) \notag \\
&= \sum_{p \in \Cal{P}_{\alpha}} | \det C_{p} |^{-1} \overline{\chi_{E}(\xi) \hat{h}_{p}(\xi)} \chi_{E} (\xi + \alpha) \hat{P g}_{p}(\xi + \alpha) \notag \\
&= \chi_{E \cap (E - \alpha)}(\xi) \sum_{p \in \Cal{P}_{\alpha}} | \det C_{p} |^{-1} \overline{\hat{h}_{p}(\xi)} \hat{g}_{p}(\xi + \alpha),
\end{align}
hence, 
\[ \sum_{p \in \Cal{P}} | \det C_{p} |^{-1} \overline{\hat{h}_{p}(\xi)} \hat{g}_{p}(\xi) = 1 \ a.e. \ \xi \in E. \]

Moreover, we must have for all $v \in V_{E}$,
\[ 0 = \sum_{p \in \Cal{P}} \sum_{k \in \Bb{Z}^d} \langle v, T_{C_p k} P g_p \rangle T_{C_p k} P^{\perp} h_p \]
which is equivalent to
\[ 0 = \sum_{p \in \Cal{P}} \sum_{k \in \Bb{Z}^d} \langle \cdot , T_{C_p k} P g_p \rangle T_{C_p k} P^{\perp} h_p. \]
Therefore, by Corollary \ref{C:Trans1}, we must have for every $\alpha \in \Lambda$ and almost every $\xi$:
\begin{align} \label{E:local2}
0 &= \sum_{p \in \Cal{P}_{\alpha}} | \det C_{p} |^{-1} \overline{\widehat{P^{\perp} h}_{p}(\xi)} \widehat{P g}_{p}(\xi + \alpha) \notag \\
&= \sum_{p \in \Cal{P}_{\alpha}} | \det C_{p} |^{-1} \overline{\chi_{\widetilde{E}}(\xi) \hat{h}_{p}(\xi)} \chi_{E} (\xi + \alpha) \hat{P g}_{p}(\xi + \alpha) \notag \\
&= \chi_{\widetilde{E} \cap (E - \alpha)}(\xi) \sum_{p \in \Cal{P}_{\alpha}} | \det C_{p} |^{-1} \overline{\hat{h}_{p}(\xi)} \hat{g}_{p}(\xi + \alpha).
\end{align}
By combining equations (\ref{E:local1}) and (\ref{E:local2}), we have
\[ \sum_{p \in \Cal{P}_{\alpha}} | \det C_{p} |^{-1} \overline{\hat{h}_{p}(\xi)} \hat{g}_{p}(\xi + \alpha) = 0 \ a.e. \ \xi \in E - \alpha.\]
\end{proof}

The following example shows that it is possible for $\{x_j\}$ to be an $M$-subspace dual to $\{y_j\}$, while $\{y_j\}$ is not an $M$-subspace dual to $\{x_j\}$.  It also shows that in the case of $M = V_{E}$ for some $E$, it is not necessary for $\Theta^{*}_{X} \Theta_{Y}$ to be in the von Neumann algebra $\{T_{\beta}:\beta \in \Bb{R}^d \}'$.

\begin{example} \label{E:bidual}
Let $M = V_{[-1/4,1/4]}$ and let $\psi$ be such that $\hat{\psi}$ is supported on $[-1/2,1/2]$, bounded, and identically $1$ on $[-1/4,1/4]$.  Define $\phi$ by $\hat{\phi}(\cdot ) = \hat{\psi}(\cdot ) + \hat{\psi}(\cdot - 1)$.  Then $\{ T_{k} \psi : k \in \Bb{Z}\}$ is $V_{[-1/4,1/4]}$-dual to $\{ T_{k} \phi : k \in \Bb{Z}\}$ but $\{ T_{k} \phi : k \in \Bb{Z}\}$ is NOT a $V_{[-1/4,1/4]}$-dual to $\{ T_{k} \psi : k \in \Bb{Z}\}$.

To see why this is the case, notice that $\overline{\hat{\phi}(\xi)}\hat{\psi}(\xi) = 1$ on $[-1/4,1/4]$.  For $k \neq 0$, $\overline{\hat{\psi}(\xi)}\hat{\phi}(\xi + k) = 0$ for $\xi \in  [-1/4,1/4] - k$.  However, for $k = - 1$, $\overline{\hat{\phi}(\xi)}\hat{\psi}(\xi - 1) \neq 0$ for $\xi \in  [-1/4,1/4] + 1$.

Alternatively, notice that for any $(c_k) \in l^2(\Bb{Z})$, $\sum_{k \in \Bb{Z}} c_k T_{k} \phi$ has Fourier transform which is repeated twice, once on $[-1/2,1/2]$ and once on$[1/2,3/2]$, whence, $\{T_{k} \phi\}$ cannot be a $V_{[-1/4,1/4]}$-dual to $\{ T_{k} \psi\}$.
\end{example}

\begin{corollary}
Let $\{T_{C_p k} g_p\}$ be Bessel and satisfy the local integrability condition, and let $E \subset \ltwod$ be measurable.  Then $\{T_{C_p k} g_p\}$ is a $V_E$-Plancherel frame if and only if 
\begin{enumerate}
\item $\displaystyle{\sum_{p \in \Cal{P}} | \det C_{p} |^{-1} |\hat{g}_{p}(\xi)|^2 = 1 \ a.e. \ \xi \in E;}$
\item for all $\alpha \in \Lambda \setminus \{0\}$,
$\displaystyle{\sum_{p \in \Cal{P}_{\alpha}} | \det C_{p} |^{-1} \overline{\hat{g}_{p}(\xi - \alpha)} \hat{g}_{p}(\xi) = 0 \ a.e. \ \xi \in E.}$
\end{enumerate}
\end{corollary}

\begin{corollary}
Suppose $A$ is an expansive integer matrix and $\Cal{U}_{A}(\Psi)$ and $\Cal{U}_{A}(\Phi)$ are Bessel sequences, and let $E$ be measurable.  Then $\Cal{U}_{A}(\Phi)$ is a $V_{E}$-subspace dual to $\Cal{U}_{A}(\Psi)$ if and only if
\begin{enumerate}
\item $\displaystyle{\sum_{i=1}^{r} \sum_{j = -\infty}^{\infty} \overline{\hat{\phi}_{i}(A^{*j}\xi)} \hat{\psi}_{i}(A^{*j}\xi) = 1 \ a.e. \ \xi \in E;}$
\item for every $q \in \Bb{Z}^d \setminus A^{*} \Bb{Z}^d$,
$\displaystyle{\sum_{i=1}^{r} \sum_{j = 0}^{\infty} \overline{\hat{\phi}_{i}(A^{*j}(\xi - q))} \hat{\psi}_{i}(A^{*j}\xi) = 0 \ a.e. \ \xi \in E}$.
\end{enumerate}
\end{corollary}
\begin{proof}
See the computation in Lemma \ref{L:tq}.
\end{proof}

\section*{Conclusion}

We have demonstrated characterization theorems for orthogonal frames consisting of regular translates, in particular affine and quasi-affine frames.  Our techniques here work in fairly general settings, including the case of Weyl-Heisenberg frames.  We have not included those results here, however, since stronger results appear in \cite{B99a,BL02a}.  Moreover, the techniques fall short with regular translation systems with different parameters and also do not apply to irregular systems.  We end the paper with a few open questions.
\begin{enumerate}
\item If the frames $\Cal{U}_{A,X}(\Psi)$ and $\Cal{U}_{B,Y}(\Phi)$ are dual, is it necessary that $X = Y$?
\item What is a full characterization of the orthogonality of $\Cal{U}_{A,X}(\Psi)$ and $\Cal{U}_{B,Y}(\Phi)$?
\item What about the case of irregular wavelet frames?
\end{enumerate}
\bibliographystyle{amsplain}
\bibliography{biblio}
\nocite{LTW03a}

\end{document}